\documentclass{IEEEtran4PSCC}

\ifCLASSINFOpdf
   \usepackage[pdftex]{graphicx}
\else
   \usepackage[dvips]{graphicx}
\fi

\usepackage[cmex10]{amsmath}

\makeatletter
\let\old@ps@headings\ps@headings
\let\old@ps@IEEEtitlepagestyle\ps@IEEEtitlepagestyle
\def\psccfooter#1{%
    \def\ps@headings{%
        \old@ps@headings%
        \def\@oddfoot{\strut\hfill#1\hfill\strut}%
        \def\@evenfoot{\strut\hfill#1\hfill\strut}%
    }%
    \def\ps@IEEEtitlepagestyle{%
        \old@ps@IEEEtitlepagestyle%
        \def\@oddfoot{\strut\hfill#1\hfill\strut}%
        \def\@evenfoot{\strut\hfill#1\hfill\strut}%
    }%
    \ps@headings%
}
\makeatother

\psccfooter{%
        \parbox{\textwidth}{\hrulefill \\ \small{23rd Power Systems Computation Conference} \hfill \begin{minipage}{0.2\textwidth}\centering \vspace*{4pt} \includegraphics[scale=0.06]{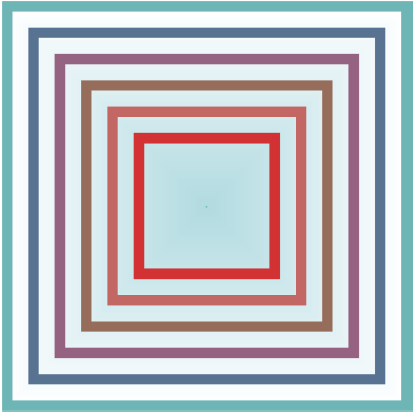}\\\small{PSCC 2024} \end{minipage} \hfill \small{Paris, France --- 2024}}%
}

\usepackage{color,soul}
\usepackage{units}
\usepackage[shortlabels]{enumitem}
\usepackage{cite}
\usepackage{lipsum}
\usepackage{booktabs}

\newcommand\blfootnote[1]{%
  \begingroup
  \renewcommand\thefootnote{}\footnote{#1}%
  \addtocounter{footnote}{-1}%
  \endgroup
}

\usepackage{amsmath,amsfonts,amsthm,amssymb}
\usepackage[mathscr]{euscript}
\usepackage[bb=boondox]{mathalfa}
\usepackage{algorithm}
\usepackage[noend]{algpseudocode}
\usepackage{array}
\usepackage{mdframed}
\usepackage{bm}

\usepackage[noabbrev,capitalize]{cleveref}

\crefname{equation}{}{}
\Crefname{equation}{}{}
\crefname{section}{}{}
\Crefname{section}{}{}

\theoremstyle{definition} 

\theoremstyle{plain} 

\theoremstyle{remark} 
\newtheorem{remark}{Remark}

\usepackage{tikz}

\newcommand{\set}[1]{\mathscr{#1}} 
\newcommand{\norm}[1]{\left\lVert#1\right\rVert} 

\newcommand*{\tran}{^{\hspace{-0.1em}\mkern-1.mu\mathsf{T}}}

\DeclareMathOperator{\cvar}{CVaR}
\DeclareMathOperator{\var}{VaR}

\DeclareMathOperator*{\argmin}{arg\,min}

\makeatletter
\newcommand{\pushright}[1]{\ifmeasuring@#1\else\omit\hfill$\displaystyle#1$\fi\ignorespaces}
\newcommand{\pushleft}[1]{\ifmeasuring@#1\else\omit$\displaystyle#1$\hfill\fi\ignorespaces}
\makeatother

\newcolumntype{C}[1]{>{\centering\arraybackslash}p{#1}} 

\newcounter{box}

\usepackage{placeins}
\newcommand{\subparagraph}{}
\usepackage{titlesec}

\begin{document}
\bstctlcite{IEEE:BSTcontrol} 

\title{Prescribed Robustness in Optimal Power Flow}


\author{
\IEEEauthorblockN{Robert Mieth, H. Vincent Poor}
\IEEEauthorblockA{Department of Electrical and Computer Engineering,
Princeton University,
Princeton, USA.
}
}


\maketitle

\begin{abstract}
For a timely decarbonization of our economy, power systems need to accommodate increasing numbers of clean but stochastic resources.
This requires new operational methods that internalize this stochasticity to ensure safety and efficiency.
This paper proposes a novel approach to compute adaptive safety intervals for each stochastic resource that internalize power flow physics and optimize the expected cost of system operations, making them ``prescriptive''.
The resulting intervals are interpretable and can be used in a tractable robust optimal power flow problem as uncertainty sets.
We use stochastic gradient descent with differentiable optimization layers to compute a mapping that obtains these intervals from a given vector of context parameters that captures the expected system state.
We demonstrate and discuss the proposed approach on two case studies.
\end{abstract}

\blfootnote{
This work was supported in part by a grant from the U.S. National Science Foundation under grant number ECCS-2039716 and a grant from the C3.ai Digital Transformation Institute. 
The work of R. Mieth was supported by a fellowship from the German Academy of Sciences, Leopoldina.
}

\section{Introduction}

The ongoing deployment of clean but stochastic energy resources challenges established power system operations by threatening their ability to ensure system security and economic efficiency. 
To tackle this, numerous effective stochastic optimizations methods have been proposed, covering, for example, system scheduling and control \cite{roald2017chance,lee2021robust,bienstock2014chance} and electricity markets and the unit commitment process \cite{kazempour2018stochastic,bertsimas2012adaptive,liang2022inertia}.
However, the adoption of \cite{roald2017chance,lee2021robust,bienstock2014chance,kazempour2018stochastic,bertsimas2012adaptive,liang2022inertia} and similar proposals is obstructed by the necessity to significantly alter established operations, an expensive and risky process for system operators.
Moreover, stochastic approaches to power system operations typically optimize reserve allocations and control policies implicitly based on probabilistic models and risk parameters defined by the system operator.
This complicates a transparent communication towards power system stakeholders (e.g., owners of generation assets, electricity market participants who are interested in anticipating price formation and scheduling procedures).
To overcome these barriers and further facilitate the adoption of clean energy technology, actionable methods for power system operation and planning that smartly internalize resource stochasticity while remaining interpretable and largely compatible with established processes are required.

Motivated by this requirement, this paper proposes a data-driven approach inspired by \cite{wang2023learning} to compute uncertainty sets for stochastic resources that (i) robustify generation and reserve allocation, (ii) internalize system physics, and (iii) minimize risk-adjusted cost while avoiding the need for a stochastic decision-making process.
Extending \cite{wang2023learning}, we introduce an additional prescription step to dynamically adjust the uncertainty sets for a given problem context. We show that this leads to significant performance increases.

\subsection{General problem formulation}

Consider power system operations as a two-stage optimization problem with uncertain parameters $\bm{\xi}$, e.g., real-time nodal demand or renewable energy injection.
Given a current parametrization of the system $\bm{\zeta}$ (e.g., net-demand forecast, {generator availability, grid topology information}), the first stage computes a resource allocation $\bm{x}$ (e.g., generation schedules, reserve allocations).  
The second stage then decides on recourse actions 
$\bm{y}$ depending on the first-stage allocation $\bm{x}^*$ and the realization of $\bm{\xi}$.
In current industry practice, the two stages are solved deterministically and independently using fixed security parameters $\bm{\theta}$ (e.g., reserve requirements) in the first stage: 
\begin{alignat}{3}
 &\text{1st stage:} \quad  \hspace{1.8em} \mathcal{P}(\bm{\zeta},\bm{\theta}) &&= \min_{\bm{x}\in\set{X}(\bm{\zeta},\bm{\theta})} \mathcal{C}^{\rm F}(\bm{x}, \bm{\zeta}, \bm{\theta}) \label{eq:first_stage}\\
 &\text{2nd stage:} \quad 
    \mathcal{Q}(\bm{x}^*, \bm{\zeta}, \bm{\xi}) &&= \min_{\bm{y}\in\set{Y}(\bm{x}^*, \bm{\zeta}, \bm{\xi})} \mathcal{C}^{\rm S} (\bm{x}^*, \bm{\zeta}, \bm{\xi})
\end{alignat}
While this process is computationally efficient and transparent, it ignores any probabilistic information that my be available from models or historic observations of $\bm{\xi}$.

As an alternative to replacing this deterministic two-step procedure with probabilistic optimization (see discussion above), this paper proposes an adaptive computation of parameters $\bm{\theta}$ that internalizes information on the stochasticity of $\bm{\xi}$ and is aware of $\bm{\zeta}$. 
We highlight that vector $\bm{\zeta}$ can be considered richer than just containing problem parameters, but may also contain additional covariates of $\bm{\xi}$ (e.g., weather information, forecasts). 
We therefore call $\bm{\zeta}$ \textit{context} parameter.
Given $\bm{\zeta}$, the optimal choice of $\bm{\theta}$ in the first stage is 
\begin{equation}
    \bm{\theta}^*(\bm{\zeta}) \in \argmin_{\bm{\theta}}\mathbb{E}_{(\bm{\xi}\mid\bm{\zeta})}\big[\mathcal{P}(\bm{\zeta},\bm{\theta}) + \mathcal{Q}(\bm{x}^*, \bm{\zeta}, \bm{\xi})\big],
\label{eq:general_problem}
\end{equation}
Following \cite{bertsimas2020predictive}, we call $\bm{\theta}^*$ \textit{prescriptive} in the context of $\bm{\zeta}$, as it provides a parametrization that minimizes the conditional expectation over $\bm{\xi}$.
To avoid re-solving \cref{eq:general_problem} every time (which might be hard or impossible), the objective of this paper is to obtain a mapping $\mathcal{M}$ that computes $\bm{\theta}^*$ from $\bm{\zeta}$, i.e.:
\begin{equation}
\min_{\mathcal{M}}\mathbb{E}_{(\bm{\zeta,\bm{\xi}})}\big[\mathcal{P}(\bm{\zeta}, \mathcal{M}(\bm{\zeta})) + \mathcal{Q}(\bm{x^*}, \bm{\zeta}, \bm{\xi})].
\label{eq:general_training_problem}
\end{equation}

Solving problem \cref{eq:general_training_problem} is hard in general and its tractability depends on the definition of the operational problem $\mathcal{P}$, $\mathcal{Q}$ and the mapping $\mathcal{M}$. 
This paper studies the following approach:
\begin{itemize}[left=0pt]
\itemsep0em
    \item For $\mathcal{P}$, $\mathcal{Q}$ we study a DC optimal power flow (OPF) problem where the first stage computes a reserve allocation alongside a second-stage control policy that ensures system safety for predefined uncertainty regions given by $\bm{\theta}$. 
    {The resulting problem is a robust optimization problem for which we compute uncertainty sets that avoid being overly conservative through a decision-aware tuning of $\bm{\theta}$. We deliberately chose to study a DC OPF with affine recourse in this paper due to its popularity and wide range of applications \cite{bienstock2014chance,roald2017chance,lubin2015robust,weinhold2023uncertainty}, allowing us to focus on the derivation and analysis of the prescriptive and decision-dependent uncertainty set tuning.}
    Section~\ref{sec:operation_model} explains the model in detail.

    \item  For $\mathcal{M}$ we focus on an interpretable linear model $\mathcal{M}_{\bm{w}}(\bm{\zeta}) = \bm{M} \bm{\zeta} + \bm{m}$ parametrized by $\bm{w} = (\bm{M}, \bm{m})$ similar to the approach in \cite{morales2023prescribing}. 
    We will refer to $\bm{w}$ as \textit{weights} throughout this paper to avoid confusion with \textit{parameters} $\bm{\zeta}$ and $\bm{\theta}$.

    \item The resulting problem is a bi-level program that we solve using stochastic gradient descent with differentiable optimization layers. We discuss this in Section~\ref{sec:solution_approach}.
\end{itemize}

\subsection{Related literature}

Robust optimization for OPF has been studied alongside numerous proposals for general stochastic OPF. Mainly introduced by the seminal work in \cite{bertsimas2012adaptive}, robust OPF can provide security with manageable computational complexity and transparent communication of the considered uncertain region for each resource. More recent variants include data-driven approaches \cite{bertsimas2018data} and tractable extensions to AC-OPF \cite{louca2018robust,louca2016stochastic,lee2021robust}. 
However, defining uncertainty sets that ensure high-quality decisions remains tricky \cite{golestaneh2018polyhedral}.
Adaptive robust programming approaches, e.g., \cite{bertsimas2012adaptive,lorca2016multistage,isuru2020piecewise}, allow constraint violation estimates, but cannot internalize the impact of the set definition on the decision outcome and cost.
Popular chance-constraint approaches, e.g., \cite{bienstock2014chance}, that rely on pre-defined ellipsoidal sets have the same weakness.

Adoption barriers for stochastic optimization in power systems have been previously highlighted in \cite{morales2014electricity,vdvorkin2018setting} and options for approximating their performance with no or little alteration of current industry practice have been explored in \cite{wang2013flexiramp,morales2014electricity,vdvorkin2018setting,garcia2021application,morales2023prescribing}. 
While \cite{wang2013flexiramp,vdvorkin2018setting} study more flexible reserve products and reserve requirements informed by an auxiliary stochastic program, respectively, \cite{morales2014electricity,garcia2021application,morales2023prescribing} propose methods that prescribe intentionally biased input parameters (forecasts) to the first stage problem. This approach fits our formulations in \cref{eq:general_problem} and \cref{eq:general_training_problem}. In \cite{morales2014electricity} the authors solve a bi-level program for each instance of the first stage to obtain a (prescriptive) alternative wind power forecast. Building on this idea, \cite{morales2023prescribing} computes a mapping from the original to a prescribed net-demand forecast. To this end, the authors include the first stage optimality conditions in the second stage to compute the map using in a single stochastic program.
Pursuing the similar objective of obtaining optimally biased forecasts that reflect the asymmetric power system cost structure (generation excess can typically be handled more cheaply than shortage), \cite{garcia2021application} propose a bi-level program with a scaleable solution heuristic. 

Models (e.g., for forecasting) that minimize the loss of a downstream optimization task have gained general popularity as \textit{end-to-end learning} \cite{kotary2021end} or \textit{smart predict-and-optimize} \cite{elmachtoub2022smart}.
Many exciting results in this direction have been unlocked by differentiable optimization frameworks, e.g., \cite{agrawal2019differentiable}, that enable efficient iterative model training procedures that internalize optimization layers.
For example, \cite{donti2017task} train a demand forecast model that minimizes the expected cost of generation excess and shortage, \cite{liang2022operation} train a generative network to obtain adversarial forecast scenarios to improve reserve allocation, \cite{wahdany2023more} create a wind power forecast model that minimizes wind spillage, and \cite{vdvorkin2023price} tune their wind power prediction model to minimize forecast errors in resulting electricity prices.

\section{Operation model}
\label{sec:operation_model}

\begin{table}[t]
\caption{Nomenclature}
\footnotesize
\centering
\begin{tabular}{c p{5.5cm}}
\toprule
\multicolumn{2}{l}{\textit{Operation model:}}\\
$D$ & Number of uncertain parameters \\
$K$ & Number of constraints \\
$G$ & Number of generators \\
$V$ & Number of nodes (vertices) in the network \\
$\bm{A}$ & Matrix of balancing control variables\\
$\bm{B}$ & Linear maps from net-injections to line flows\\
$\bm{c}$ & Vector of cost parameters \\
$\bm{d}$ & Vector of active power demand \\
$\bm{f}$ & Vector of power flows \\
$\bm{p}$ & Vector of active power generation\\
$\bm{r}^+$/$\bm{r}^-$ & Vector of upward/downward generation reserves \\
$\bm{u}$ & Vector of forecasts of uncertain injections \\
$\bm{x}$ & Vector collecting model decision variables\\
$\bm{\xi}$  & Vector of forecast errors  \\

\multicolumn{2}{l}{\textit{Other variables and parameters:}} \\
$G$ & Subgradient of loss function \\
$H$ & Empirical CVaR \\
$L$ & Loss function \\
$\bm{M}$ & Affine part of linear prescription map \\
$\mathcal{M}_{\bm{w}}$ & Prescription map parameterized by $\bm{w}$ \\
$N$ & Number of samples \\
$\bm{m}$ & Constant part of linear prescription map \\
$\bm{w}$ & Parameters (weights) of prescription map \\
$\bm{\Sigma}_{\bm{u}}$ & Covariance matrix of $\bm{\xi}$ conditional to $\bm{u}$ \\
$\gamma$ & Target probability of exceedance \\ 
$\bm{\zeta}$  & Vector of context parameters \\
$\kappa$ & Step size for Lagrangian in CVaR loss \\
$\lambda$ & Lagrangian multiplier for CVaR loss \\
$\bm{\mu}$ & Center of box uncertainty set \\
$\bm{\rho}$ & Learning rate \\ 
$\bm{\sigma}$ & Vector collecting width for each dimension of the box uncertainty set \\
$\tau$ & Auxiliary variable for CVaR \\
$\bm{\theta}$ & Uncertainty set parameters \\

\multicolumn{2}{l}{\textit{Other notation:}} \\
$\cdot^v$/$\cdot^z$ & Superscript indicating iteration epoch/batch \\
$[\cdot]^+$ & $\max\{0,\cdot\}$ \\
$\bm{1}_X$ & Vector of ones of dimension $X$ \\
\bottomrule
\end{tabular}
\end{table}

We consider a short-term generator dispatch problem with balancing control (e.g., automatic generator control) with uncertain injections from stochastic wind generators.
We model these uncertain injections as a $D$-dimensional vector $\bm{u}(\bm{\xi})=\bm{u} + \bm{\xi}$, where $\bm{u}$ is a (deterministic) forecast and $\bm{\xi}$ is a vector of random forecast errors. 
Power imbalance caused by forecast errors $\bm{\xi}$ is corrected by controllable generators.
The generator output is a $G$-dimensional vector $\bm{p}(\bm{\xi}) = \bm{p} - \bm{A}\bm{\xi}$, {where $\bm{A}$ is a matrix of decision variables that defines an affine balancing control, i.e., it defines how each generator adjusts its output as a reaction to imbalance $\bm{\xi}$.}
Ensuring $\bm{A}\tran\bm{1}_{G} = \bm{1}_D$, where $\bm{1}_{D}$ is a $D$-dimensional vector of ones, ensures system balance.

\subsection{Robust optimal power flow with affine recourse}
\label{ssec:robust_opf}

To immunize the system against uncertain injections $\bm{u}(\bm{\xi})$, $\bm{p}(\bm{\xi})$ and the resulting uncertain power flows, the system operator defines an uncertainty set $\set{U}$ that captures all outcomes of $\bm{\xi}$ for which all system constrains should hold.
Given a parametrization $\bm{\zeta}$  (i.e., forecast $\bm{u}$ and demand vector $\bm{d}$---we assume in this paper that cost and system topology remain constant) the system operator solves the following robust OPF problem to decide on the generator dispatch $\bm{p}$, {generator reserves $\bm{r^+}$, $\bm{r^-}$, and transmission line security margins
$\bm{f}^{\rm RAM+}, \bm{f}^{\rm RAM-}$}:
\allowdisplaybreaks
\begin{subequations}
\begin{align}
\min \quad  & (\bm{c}^{\rm E})\tran\bm{p} + (\bm{c}^{\rm R})\tran(\bm{r}^+ + \bm{r}^{-}) \label{base_dcopf:objective}\\
\text{s.t.}\quad 
    & \bm{1}_G\tran\bm{p}  = \bm{1}_V\tran\bm{d} - \bm{1}_D\tran\bm{u} \label{base_dcopf:enerbal} \\  
    & \bm{A}\tran\bm{1}_G= \bm{1}_D \label{base_dcopf:resbal}\\    
    & \bm{p} + \bm{r}^+ \le \bm{p}^{\rm max} \label{base_dcopf:gen_uplim}\\   
    & \bm{p} - \bm{r}^- \ge \bm{p}^{\rm min} \label{base_dcopf:gen_dnlim}\\
    & \bm{B}^{\rm G} \bm{p} + \bm{B}^{\rm W} \bm{u} - \bm{B}^{\rm B}\bm{d} \le \bm{f}^{\rm max} - \bm{f}^{\rm RAM+} \label{base_dcopf:lin_uplim}\\ 
    &-(\bm{B}^{\rm G} \bm{p} + \bm{B}^{\rm W} \bm{u} - \bm{B}^{\rm B}\bm{d}) \le \bm{f}^{\rm max} - \bm{f}^{\rm RAM-} \label{base_dcopf:lin_dnlim} \\ 
    & -\bm{A}\bm{\xi} \le \bm{r}^+ &&\hspace{-2cm} \forall \bm{\xi}\in\set{U} \label{base_dcopf:rob_resconst_up}\\
    & \bm{A}\bm{\xi} \le \bm{r}^-  &&\hspace{-2cm} \forall \bm{\xi}\in\set{U} \label{base_dcopf:rob_resconst_lo}\\
    & (\bm{B}^{\rm W} - \bm{B}^{\rm G}\bm{A})\bm{\xi} \le \bm{f}^{\rm RAM+}\!  &&\hspace{-2cm} \forall \bm{\xi}\in\set{U} \label{base_dcopf:rob_flowconst_up}\\
    & -(\bm{B}^{\rm W} - \bm{B}^{\rm G}\bm{A})\bm{\xi} \le \bm{f}^{\rm RAM-}\! &&\hspace{-2cm} \forall \bm{\xi}\in\set{U} \label{base_dcopf:rob_flowconst_lo}
\end{align}%
\label{eq:base_dcopf}%
\end{subequations}%
\allowdisplaybreaks[0]%
The objective~\cref{base_dcopf:objective} minimizes system cost given energy cost and reserve provision cost vectors $\bm{c}^{\rm E}$ and $\bm{c}^{\rm R}$.
Energy balance \cref{base_dcopf:enerbal} ensures that the total generator injections equals the total system demand $\bm{d}$ and $\bm{u}$.
Similarly, \cref{base_dcopf:resbal} ensures that  all forecast errors are balanced. 
Constraints \cref{base_dcopf:gen_uplim,base_dcopf:gen_dnlim} enforce the technical production limits of each controllable generator. 
Constraints \cref{base_dcopf:lin_uplim,base_dcopf:lin_dnlim} map the power injections and withdrawals of each resource and load to a resulting power flow {via linear maps $B^{\rm G}$, $B^{\rm W}$, $B^{\rm B}$, e.g., obtained from the DC power flow approximation \cite{bolognani2015fast}}. 
Vectors $\bm{f}^{\rm RAM+}$ and $\bm{f}^{\rm RAM-}$ are the remaining available margins for each power transmission line, i.e., the difference between the power flow caused by the forecast injections and the upper and lower line limits.
Constraints \cref{base_dcopf:rob_resconst_lo,base_dcopf:rob_resconst_up,base_dcopf:rob_flowconst_lo,base_dcopf:rob_flowconst_up} enforce robust constraints on the system response to uncertain forecast errors $\bm{\xi}$.
Constraints \cref{base_dcopf:rob_resconst_lo,base_dcopf:rob_resconst_up} ensure that generator balancing responses do not exceed the available reserves and \cref{base_dcopf:rob_flowconst_lo,base_dcopf:rob_flowconst_up} ensure that the resulting power flow changes do not exceed the remaining available margins on each power line for any $\bm{\xi}\in\set{U}$.

Problem \cref{eq:base_dcopf} can be re-written in a more concise form. 
We collect all decision variables in a vector $\bm{x}$, cost vectors $\bm{c}^{\rm E}$, $\bm{c}^{\rm R}$ in a vector $\bm{c}$, denote the feasible space defined by constraints \cref{base_dcopf:resbal,base_dcopf:enerbal,base_dcopf:gen_uplim,base_dcopf:gen_dnlim,base_dcopf:lin_uplim,base_dcopf:lin_dnlim} as $\set{F}(\bm{\zeta})$ and write: 
\begin{subequations}
\begin{align}
\min_{\bm{x}\in\set{F}(\bm{\zeta})}\ 
    &  \bm{c}\tran\bm{x}  \label{compact_dcopf:objective} \\
\text{s.t.}\quad
   & \max_{k=1,...,K}[ \bm{a}_{k}\tran\bm{\xi} + b_k] \le 0, \quad \forall\bm{\xi}\in\set{U} \label{compact_dcopf:max_of_affine_robust}
\end{align}%
\label{eq:compact_dcopf}%
\end{subequations}%
where $\bm{a}_k$ and $b_k$ are, respectively, the $k$-th row and $k$-th entry of the 
$K\times D$ matrix and $K\times 1$ vector
\begin{equation*}
    \begin{bmatrix}
    -\bm{A} \\ \bm{A} \\ (\bm{B}^{\rm\bm{W}} - \bm{B}^{\rm G}\bm{A}) \\ -(\bm{B}^{\rm W} - \bm{B}^{\rm G}\bm{A})
    \end{bmatrix}  \text{   and   } 
    \begin{bmatrix}
    -\bm{r}^+ \\ -\bm{r}^- \\ - \bm{f}^{\rm RAM+} \\ -\bm{f}^{\rm RAM-}
    \end{bmatrix}.
\end{equation*} 
Note that \cref{compact_dcopf:max_of_affine_robust} is an exact reformulation of \cref{base_dcopf:rob_resconst_lo,base_dcopf:rob_resconst_up,base_dcopf:rob_flowconst_lo,base_dcopf:rob_flowconst_up} as the maximum of $K$ affine functions.

\subsection{Uncertainty set formulation}
Model \cref{eq:base_dcopf} cannot be solved directly but requires a definition of $\set{U}$ alongside a tractable reformulation of constraints \cref{base_dcopf:rob_resconst_lo,base_dcopf:rob_resconst_up,base_dcopf:rob_flowconst_lo,base_dcopf:rob_flowconst_up}.
This paper focuses on box uncertainty sets, as they provide a clear safety region for each uncertain resource. 
Other formulations are possible \cite{bertsimas2011theory,gorissen2013robust}.
Box uncertainty sets ensure that constraints are feasible for a security interval along each dimension of the uncertain vector $\bm{\xi}$. 
We can define such a set as $\set{U}^{\rm box}(\bm{\theta}) = \{\xi_j|\xi_j\in[\mu_j-\sigma_j, \mu_j+\sigma_j],\ j=1,...,D\}$ parametrized by $\bm{\theta}=(\bm{\mu}, \bm{\sigma})$ with $\bm{\mu}=[\mu_1,...,\mu_D]$ and $\bm{\sigma}=[\sigma_1,...,\sigma_D]$. Parameter $\bm{\mu}$ is the center of the security interval and can be interpreted as a forecast error bias. 
Parameter $\bm{\sigma}$ defines the width of the interval.
Using $\set{U} = \set{U}^{\rm box}$ and introducing auxiliary variable $\bm{t}_k$, \cref{compact_dcopf:max_of_affine_robust} becomes
\begin{equation}
\begin{aligned}
        \bm{a}_k\tran \bm{\mu} + \bm{t}_k\tran\bm{\sigma} \le b_k, \ 
        \bm{t}_k \ge |\bm{a}_k|, \quad \forall k=1,..,K. 
\label{eq:robust_reformulation}
\end{aligned}
\end{equation}

\subsection{Real-time cost and security}
\label{ssec:real_time_cost}
The choice of uncertainty set parameters $\bm{\theta}$ implies a trade-off between security in real-time and cost in the first-stage decision. 
For example, choosing $\bm{\theta}$ such that $\set{U}$ is large and covers all potential outcomes of $\bm{\xi}$ will lead to high security in real-time but also to high first-stage cost. 
A good choice of $\bm{\theta}$ will balance this trade-off by minimizing the combined first- and second-stage cost as defined in \cref{eq:general_problem}. We consider two relevant approaches to quantify security of the robust problem from Section~\ref{ssec:robust_opf} in real-time.

\subsubsection{Cost of exceedance}
We can define the real-time cost as cost of exceedance by imposing a penalty for insufficient reserves $\bm{r}^+$, $\bm{r}^-$, $\bm{f}^{\rm RAM +}$, $\bm{f}^{\rm RAM -}$. Using the notation from \cref{eq:compact_dcopf}, we compute this cost as
\begin{equation}
   \mathcal{C}^{\rm S}(\bm{x}, \bm{\xi}) = \sum_{k=1}^K c_k^{\rm viol}\big[ \bm{a}_{k}\tran\bm{\xi} + b_k\big]^+ ,
\end{equation}
where $[\cdot]^+ = \max\{\cdot, 0\}$ and $c_k^{\rm viol}$ is the cost for exceeding the reserve given by $b_k$.
Cost $c_k^{\rm viol}$ could, for example, reflect the cost of procuring emergency resources or load shedding. 

\subsubsection{Probability of exceedance}
Instead of minimizing cost of exceedance, the system operator may be interested in a probabilistic guarantee that real-time operations do not exceed reserves, i.e.:
\begin{equation}
\mathbb{P}\Big[\max_{k=1,...,K}[ \bm{a}_{k}\tran\bm{\xi} + b_k] < 0 \Big] \ge 1-\gamma.
\end{equation}
Here, $(1-\gamma)$ defines the target probability of no constraint exceeding its limits and $\gamma$ is a small risk factor.

\section{Solution Approach}
\label{sec:solution_approach}
We now show an iterative method inspired by the results in \cite{wang2023learning} to obtain the desired mapping $\mathcal{M}_{\bm{w}}(\bm{\zeta})=\bm{M}\bm{\zeta} + \bm{m}$ that returns uncertainty set parameters $\bm{\theta}$, which, in turn, optimally parametrize the first-stage uncertainty set with respect to the cost and probability of exceedance discussed in Section~\ref{ssec:real_time_cost}.

\subsection{Cost of exceedance}
\label{ssec:coe_solution_approach}

The problem to compute the optimal choice of weights $\bm{w}=(\bm{M}, \bm{m})$ that minimize the combined expected first- and second-stage cost is the bi-level problem:
\begin{subequations}
\begin{align}
    \min_{\bm{w}}\quad & \mathbb{E}_{(\bm{\zeta},\bm{\xi})}\big[\bm{c}\tran\bm{x}^* + \mathcal{C}^{\rm S}(\bm{x}^*, \bm{\xi})\big] \label{eq:coe_objective}\\
\text{s.t.}\quad 
    & \bm{\theta} = \mathcal{M}_{\bm{w}}(\bm{\zeta}) \label{eq:coe_prescription}\\
    & \hspace{-0.6cm}\bm{x}^*(\bm{\zeta},\bm{\theta})\! \in \!
\left\{\begin{array}{l}
    \!\!\argmin_{\bm{x}\in\set{F}(\bm{\zeta})}\ \!\bm{c}\tran\bm{x} \\[1.3ex]
    \!\!\text{s.t.}\ \max_{k}[ \bm{a}_{k}\tran\bm{\xi} + b_k]\! \le\! 0,\ \forall\bm{\xi}\!\in\!\set{U}(\bm{\theta}).
\end{array}
  \right. \label{eq:coe_inner}
\end{align}%
\label{eq:coe_bilevel_problem}%
\end{subequations}%
We solve \cref{eq:coe_bilevel_problem} using a stochastic gradient descent approach (see, e.g., \cite{amari1993backpropagation}) with $v^{\rm max}$ outer steps (epochs) indexed by $v$ and $z^{\rm max}$ {inner steps (``mini-batches'' \cite{gower2019sgd})} indexed by $z$.
For each step $z$ we assume having access to an \textit{individual} context parameter sample $\bm{\zeta}^z$ (e.g., obtained from historic observations or a sample generation mechanism) and a \textit{set of $N_z$ samples} $\set{X}^z = \{\bm{\xi}_{i}^z\}_{i=1}^{N_z}$ of $\bm{\xi}$ (again, either obtained from historical observations or a sample generation mechanism). 
We note that $\set{X}^z$ may be conditional to $\bm{\zeta}^z$. See, for example, \cite{dvorkin2015uncertainty} who also provide an approach to sample $\set{X}^z$ for a given wind power forecast $\bm{u}^z$ (which is part of the parametrization $\bm{\zeta}^z$).
We provide additional discussion on the relationship between $\bm{\zeta}^z$ and $\set{X}^z$ in the case study below.
We define $\bm{x}^z = \bm{x}^*(\bm{\zeta}^z, \bm{\theta}^z)$ as the solution to the inner problem \cref{eq:coe_inner}. The loss function $L(\bm{w}^z; \bm{\zeta}^z, \bm{\xi}^z)$ corresponds to the problem objective in \cref{eq:coe_inner} and we compute it as the empirical mean given $\set{X}^z$:
\begin{equation}
    L^{\rm C}\!(\bm{w}^z\!; \bm{\zeta}^z\!, \set{X}^z\!)\! =\! \bm{c}\tran\bm{x}^z\! +\! \frac{1}{N_z}\! \sum_{i=1}^{N_z}  \sum_{k=1}^K\! c_k\big[ (\bm{a}_{k}^z)\tran\bm{\xi}_{i}^z\! +\! b_k^z\big]^+.
\end{equation}
{Recall from \cref{compact_dcopf:max_of_affine_robust} that $\bm{a}_k^z$ and $b_k^z$ are part of $\bm{x}^z$.}  
Following \cite{wang2023learning} we estimate the derivative of $L^{\rm C}(\bm{w}^z; \bm{\zeta}^z, \bm{\xi}^z)$ using a subgradient $G^{\rm C}(\bm{w}^z; \bm{\zeta}^z, \set{X}^z)$ computed over samples of $\bm{\zeta}$ and $\bm{\xi}$.
Notably, the computation of $G^{\rm C}(\bm{w}^z; \bm{\zeta}^z, \set{X}^z)$ requires computing gradient $\nabla_{\bm{\theta}}\bm{x}^*(\bm{\zeta}, \bm{\theta})$, i.e., the derivative of the decision variables of the inner problem over the uncertainty set parameters. We discuss an efficient approach to obtain this gradient in Section~\ref{ssec:implementation} below.
Choosing initial weights $(\bm{M}^{\rm init}, \bm{m}^{\rm init})$ and a learning rate $\rho$ the resulting solution steps are itemized in Algorithm~\ref{alg:coe_sgd}.
\begin{algorithm}
\caption{Stochastic gradient descent to solve \cref{eq:coe_bilevel_problem}}
\label{alg:coe_sgd}
\begin{algorithmic}[1]
    \State {\bf given}  $\bm{w}^0 = (\bm{M}^{\rm init}, \bm{m}^{\rm init})$,  learning rate $\rho$ 
    \For{$v=1,...,v^{\rm max}$}  \Comment{\textit{outer epoch}}
        \For{$z=1,...,z^{\rm max}$}  \Comment{\textit{inner mini-batch}}
            \State $\bm{\zeta}^z \gets$ sample context parameter
            \State $\set{X}^z \gets$ sample forecast errors conditional to $\bm{\zeta}^z$
            \State $\bm{\theta}^z \gets \mathcal{M}_{\bm{w}^{v-1}}(\bm{\zeta}^z)$ \Comment{\textit{prescribe parameters}}
            \State $\bm{x}^z \gets$ solve inner problem \cref{eq:coe_inner}
            \State $\bm{g}^z \gets G^{\rm C}(\bm{w}^{v-1}; \bm{\zeta}^z, \set{X}^z)$           
        \EndFor
        \State $\bm{g}^v \gets \frac{1}{z^{\rm max}}\sum_{z=1}^{z^{\rm max}}\bm{g}^z$ \Comment{\textit{update gradient}}
        \State $\bm{w}^{v} \gets \bm{w}^{v-1}\! -\! \rho \bm{g}^v$ \Comment{\textit{update variables}}
    \EndFor
    \State {\bf return} $\bm{w}^{v^{\rm max}}$
\end{algorithmic}
\end{algorithm}

\subsection{Probability of exceedance}
\label{ssec:poe_solution_approach}

The problem to compute the optimal choice of weights $\bm{w}$ that ensure real-time operations do not exceed their limits with a probability of at least $(1-\gamma)$ is the bi-level problem:
\allowdisplaybreaks
\begin{subequations}
\begin{align}
    \min_{\bm{w}}\quad & \mathbb{E}_{(\bm{\zeta},\bm{\xi})}\big[\bm{c}\tran\bm{x}^* \big] \\
\text{s.t.}\quad 
    &  \mathbb{P}\Big[\max_{k=1,...,K}[ (\bm{a}_{k}^*)\tran\bm{\xi} + b_k^*] < 0 \Big] \ge 1-\gamma \label{eq:poe_chance_constraint} \\
    & \big[\text{$\bm{\theta}$ and $\bm{x}^*(\bm{\zeta, \bm{\theta}})$ as in \cref{eq:coe_prescription,eq:coe_inner}}\big].
\end{align}%
\label{eq:poe_bilevel_problem}%
\end{subequations}%
\allowdisplaybreaks[0]%
Constraint \cref{eq:poe_chance_constraint} is generally non-convex and complicates the solution of \cref{eq:poe_bilevel_problem}.
To create a tractable problem, we reformulate \cref{eq:poe_chance_constraint} using conditional value-at-risk \cite{wang2023learning,mieth2023data}:
\begin{equation}
\begin{aligned}
    & \mathbb{P}\big[\max_{k=1,...,K}[ \bm{a}_{k}\tran\bm{\xi} + b_k \le 0 \big] \ge 1-\gamma \\
\Leftrightarrow \quad
    &\var_{\gamma}\big( \max_{k=1,...,K}[ \bm{a}_{k}\tran\bm{\xi} + b_k] \big) \le 0 \\
\Leftarrow \quad 
    & \cvar_{\gamma}\big( \max_{k=1,...,K}[ \bm{a}_{k}\tran\bm{\xi} + b_k ]\big) \le 0 \\
\end{aligned}
\label{eq:cvar_reformulation}
\end{equation}
where $\var_{\gamma}$ and $\cvar_{\gamma}$ are the value-at-risk and conditional value-at-risk at risk level $\gamma$. Reformulation \cref{eq:cvar_reformulation} utilizes the fact that limiting CVaR implies a limit on VaR \cite{rockafellar2000optimization}. CVaR further allows the convex reformulation \cite{rockafellar2000optimization}
\begin{equation}
    \cvar_{\gamma} = \inf_{\tau}\Big\{\mathbb{E}\Big[\frac{1}{\gamma}\big[ \max_{k=1,...,K}[ \bm{a}_{k}\tran\bm{\xi} + b_k] - \tau\big]^+ + \tau\Big]\Big\},
\label{eq:cvar_as_inf}
\end{equation}
which we use to define 
\begin{equation}
    h(\bm{x}^*, \tau, \bm{\xi}) = \frac{1}{\gamma}\big[  \max_{k=1,...,K}[ (\bm{a}^*_{k})\tran\bm{\xi} + b^*_k] - \tau\big]^+ + \tau,
\label{eq:h_function}
\end{equation}
and to reformulate \cref{eq:poe_bilevel_problem} as \cite{wang2023learning}:
\allowdisplaybreaks
\begin{subequations}
\begin{align}
\min_{\bm{w}, \tau}\quad 
    & \mathbb{E}_{(\bm{\zeta}, \bm{\xi})} \big[\bm{c}\tran \bm{x}^* \big] \\
\text{s.t.}\quad 
    & \mathbb{E}_{(\bm{\zeta}, \bm{\xi})} \big[h(\bm{x}^*, \tau, \bm{\xi}) \big] = 0 \label{eq:poe_reform_cvar} \\
    & \big[\text{$\bm{\theta}$ and $\bm{x}^*(\bm{\zeta, \bm{\theta}})$ as in \cref{eq:coe_prescription,eq:coe_inner}}\big].
\end{align}%
\label{eq:poe_reformulated}%
\end{subequations}%
\allowdisplaybreaks[0]%
We note that \cref{eq:poe_reformulated} has an additional auxiliary variable $\tau$ related to the CVaR reformulation \cref{eq:cvar_as_inf}. Also, following the logic in \cite{wang2023learning}, we note that \cref{eq:poe_reform_cvar} is an equality constraint, because zero is the optimal CVaR target given \cref{eq:cvar_reformulation} and the convexity of \cref{eq:cvar_as_inf}. 
{Finally, we note that for fixed $\bm{x}^*$ the  term $\max_{k}[ (\bm{a}^*_{k})\tran\bm{\xi} + b^*_k]$ in \cref{eq:h_function} computes the constraint violation resulting from a given $\bm{\xi}$ and may be larger than zero.}

Similar to he procedure outlined in \cite{wang2023learning,zhang2022solving} we can now solve the equality-constrained problem \cref{eq:poe_reformulated} by introducing Lagrangian multiplier $\lambda$ and defining the loss function
\begin{equation}
    L^{\rm P}((\bm{w}\!,\tau); \lambda\!,  \bm{\zeta}\!, \bm{\xi})\! = \!\mathbb{E}_{(\bm{\zeta}, \bm{\xi})} \big[\bm{c}\tran \bm{x}^* \big]\! +\! \lambda \mathbb{E}_{(\bm{\zeta}, \bm{\xi})} \big[h(\bm{x}^*\!, \tau\!, \bm{\xi}) \big]. 
\label{eq:lagrangian_expectation}
\end{equation}
Using the notation introduced in Section~\ref{ssec:coe_solution_approach} above, we compute \cref{eq:lagrangian_expectation} in each step $z$ as 
\allowdisplaybreaks
\begin{align}
    &L^{\rm P}((\bm{w}^z, \tau^z); \lambda^z, \bm{\zeta}^z, \set{X}^z) =  \label{eq:poe_loss_function}\\
    &\quad \bm{c}\tran \bm{x}^z\! +\! \lambda^z \underbrace{\frac{1}{N_z}\sum_{i=1}^{N_z}\! \! \Big( \frac{1}{\gamma}\big[\max_{k=1,...,K}[ (\bm{a}_{k}^z)\tran\bm{\xi}_i^z\! +\! b_k^z]\! -\! \tau^z\big]^+\!\! + \!\tau^z \Big)}_{H((\bm{w}^z, \tau^z); \bm{\zeta}^z, \set{X}^z)}. \nonumber
\end{align}%
\allowdisplaybreaks[0]%
Denoting $G^{\rm P}((\bm{w}^z, \tau^z); \lambda^z, \bm{\zeta}^z, \set{X}^z)$ as the gradient of \cref{eq:poe_loss_function}, $\tau^{\rm init}$, $\lambda^{\rm init}$ as initial values for $\tau$, $\lambda$, and $\kappa$ as the step size for $\lambda$ the resulting solution steps are itemized in Algorithm~\ref{alg:poe_sgd}.
\begin{algorithm}
\caption{Stochastic gradient descent to solve \cref{eq:poe_bilevel_problem}}
\label{alg:poe_sgd}
\begin{algorithmic}[1]
    \State {\bf given}  $\bm{w}^0 = (\bm{M}^{\rm init}, \bm{m}^{\rm init})$, $\tau^0=\tau^{\rm init}$, $\lambda^0 = \lambda^{\rm init}$, $\rho$, $\kappa$
    \For{$v=1,...,v^{\rm max}$}  \Comment{\textit{outer epoch}}
        \For{$z=1,...,z^{\rm max}$}  \Comment{\textit{inner mini-batch}}
            \State $\bm{\zeta}^z \gets$ sample context parameter
            \State $\set{X}^z \gets$ sample forecast errors conditional to $\bm{\zeta}^z$
            \State $\bm{\theta}^z \gets \mathcal{M}_{\bm{w}^{v-1}}(\bm{\zeta}^z)$ \Comment{\textit{prescribe parameters}}
            \State $\bm{x}^z \gets$ solve inner problem \cref{eq:coe_inner}
            \State $\bm{g}^z \gets G^{\rm C}((\bm{w}^{v-1}, \tau^{v-1}); \lambda^{v-1}, \bm{\zeta}^z, \set{X}^z)$          
            \State $H^v \gets H((\bm{w}^z, \tau^z), \bm{\zeta}^z, \set{X}^z)$ \Comment{\textit{See \cref{eq:poe_loss_function}}}
        \EndFor
        \State $\bm{g}^v \gets \frac{1}{z^{\rm max}}\sum_{z=1}^{z^{\rm max}}\bm{g}^z$ \Comment{\textit{update gradient}} 
        \State $(\bm{w}^{v}, \tau^v) \gets (\bm{w}^{v-1}, \tau^{v-1})\! -\! \rho \bm{g}^v$ \Comment{\textit{update variables}}
        \State $H^v = \frac{1}{z^{\rm max}}\sum_{z=1}^{z^{\rm max}}H^z$ \Comment{\textit{compute current CVaR}}
        \State $\lambda^{v} \gets \lambda^{v-1} + \kappa H^v$ \Comment{\textit{update Lagrangian multiplier}}
    \EndFor
    \State {\bf return} $(\bm{w}^{v^{\rm max}}, \tau^{v^{\rm max}})$
\end{algorithmic}
\end{algorithm}

{
\begin{remark}
We highlight that the mapping $\mathcal{M}_{\bm{w}}$ is only used to define the uncertainty set $\set{U}$. Once this set is obtained, problem \cref{eq:base_dcopf} can be solved by off-the-shelf solvers using reformulation \cref{eq:robust_reformulation}. This guarantees optimality and feasibility for a given $\set{U}$.
\end{remark}
}

\subsection{Implementation}
\label{ssec:implementation}

\begin{figure}
    \centering
    \includegraphics[width=0.85\linewidth]{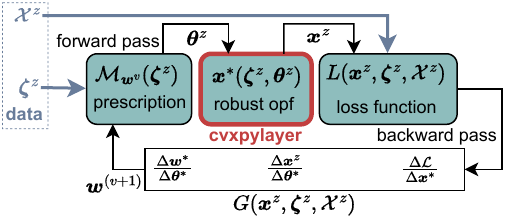}
    \caption{Computation pipeline for the steps itemized in Algorithms~\ref{alg:coe_sgd} and \ref{alg:poe_sgd}. The differentiable optimization layer using \texttt{cvxpylayers} \cite{agrawal2019differentiable}, allows the computation of subgradient $G = \nabla_{\bm{w}}L$ with a single backwards pass.}
    \label{fig:dolayers_pipeline}
\end{figure}

We implement Algorithms~\ref{alg:coe_sgd} and \ref{alg:poe_sgd} in Python using state-of-the-art machine learning packages (PyTorch) and results from \cite{agrawal2019differentiable} (package Cvxpylayers) that allow numerical differentiation through the inner optimization \cref{eq:coe_inner}.
The computation steps are illustrated in Fig.~\ref{fig:dolayers_pipeline}. The differentiable optimization layer obtains the required gradient $\nabla_{\bm{\theta}}\bm{x}^*(\bm{\zeta}, \bm{\theta})$ by differentiating through the Karush-Kuhn-Tucker conditions of the inner problem at optimality. This is performed efficiently by utilizing the implicit function theorem and by solving an inner quadratic program to obtain the required matrix inversion \cite{agrawal2019differentiating}.

We add the following implementation remarks: 
\begin{enumerate}[label=\textbf{IR{\arabic*}:}, ref=IR{\arabic*}, align=left, leftmargin=1em]
\item Formulation \cref{eq:base_dcopf} limits the decision variables $\bm{A}$ to the interval $[0,1]$. This leads to $\bm{A}$ taking small values compared to the other decision variables causing conditioning problems in the Cvxpylayers solver ECOS \cite{domahidi2013ecos}. Introducing a scaling factor for $\bm{A}$ solves this problem consistently. \label{ir:conditioniong}

\item The inner problem has to be always feasible. To ensure this, we add a slack variable to \cref{base_dcopf:enerbal} to allow the curtailment of wind power $\bm{u}$ if needed, and a slack variable to \cref{compact_dcopf:max_of_affine_robust} to account for the case that the required uncertainty set can not be met with the available reserves. Both sets of slack variables are penalized in the problem objective and conditioned using a scaling factor as in \ref{ir:conditioniong}.

\item Formulation \cref{eq:base_dcopf} is a linear problem and as such may not generally be differentiable with respect to its parameters \cite{wilder2019melding}.
This can be avoided by introducing a regularization $\pi\norm{\bm{x}}_2^2$ with regularization factor $\pi$ to the objective of the problem.
We note that in contrast to \cite{wilder2019melding}, our problem is linear with \textit{continuous} variables and adding a regularization term with a small $\pi$ did not significantly impact our results.

\end{enumerate}

Lastly, we note that we implement $\mathcal{M}_{\bm{w}}$ as two linear models $\mathcal{M}^{\bm{\mu}}=\bm{M}^{\mu}\bm{\zeta} + \bm{m}^{\mu}$ and $\mathcal{M}^{\bm{\sigma}}=\bm{M}^{\sigma}\bm{\zeta} + \bm{m}^{\sigma}$.

\section{Numerical Experiments}

\subsection{Illustrative 5-bus case}

\begin{figure}
    \centering
    \includegraphics[width=0.75\linewidth]{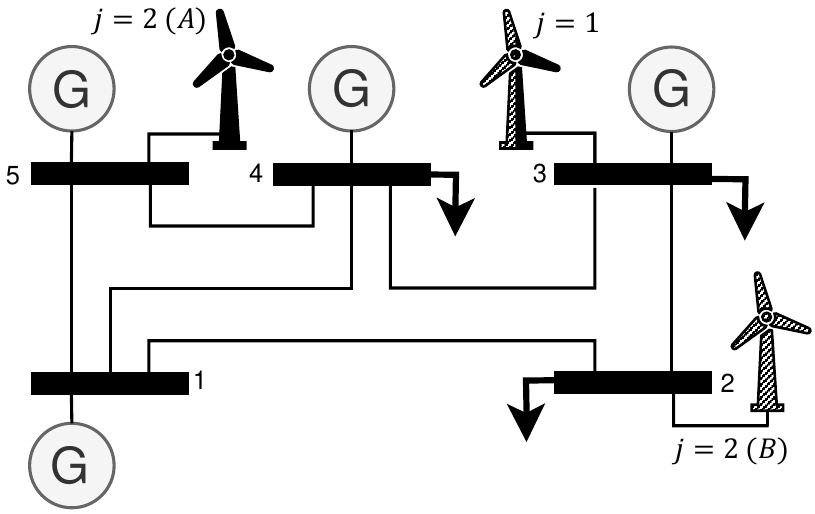}
    \caption{Schematic of the 5 bus test system in configuration A (wind farms at buses 3 and 5) and B (wind farms at buses 2 and 3).}
    \label{fig:5bus_system}
\end{figure}

We first illustrate the suggested approach using synthetic data on the small-scale ``case5'' data set from MATPOWER \cite{matpowercase5}. 
The system topology with two configurations of wind generators is shown in Fig.~\ref{fig:5bus_system}. Its parameters follow \cite{mieth2023data} with \unit[$\bm{f}^{\rm max}=(3.2, 1.52, 1.76, 0.8, 0.8, 1.92)$]{p.u.}, $\bm{c}^{\rm E}=(14,15,30,40,10)$\nicefrac{\$}{MW}, and $\bm{c}^{\rm R}=(80,80,15,30,80)$\nicefrac{\$}{MW}.

For training and testing the proposed approach we create a collection of $N=2000$ samples of $\bm{\zeta}=(\bm{d},\bm{u})$ with corresponding samples of $\bm{\xi}$ as follows.  
First, we define the nominal demand as \unit[$\bm{d}_0=(0.0, 3.0, 3.0, 4.0, 0.0)$]{p.u.} and the nominal wind forecast as \unit[$\bm{u}_0=(1.0, 1.5)$]{p.u.}. We then create $N$ samples of $\bm{d}$ by uniformly drawing from the set $\{\bm{d}\in\mathbb{R}\mid 0.5\bm{d}_0 \le \bm{d} \le 1.1\bm{d}_0\}$ and $N$ samples of $\bm{u}$ by uniformly drawing from the set $\{\bm{u}\in\mathbb{R}\mid 0.5\bm{u}_0 \le \bm{u} \le 1.1\bm{u}_0\}$.
Next, for each available sample of $\bm{\zeta}=(\bm{d},\bm{u})$ we draw a sample of forecast errors $\bm{\xi}$ from a multivariate normal distribution with a uniform correlation coefficient of $\phi=0.5$ between all wind farms. Following the observations in \cite{dvorkin2015uncertainty}, we model the forecast error standard deviation to be conditional to the forecast and the forecast error mean to be zero. The resulting normal distribution is given as:
\begin{equation}
\setlength\arraycolsep{2pt}
    \bm{\xi}\,|\,\bm{u} \!\sim\! \mathcal{N}(\bm{0}, \bm{\Sigma}_{\bm{u}}),\ \bm{\Sigma}_{\bm{u}} \!= \!\!
    \begin{bmatrix}(0.15u_1)^2  & 0.15^2\phi u_1 u_2) \\ 0.15^2\phi u_1 u_2) & (0.15u_2)^2 \end{bmatrix}\!.
\label{eq:forecast_error_sample_distribution}
\end{equation}
Finally, we set the maximum wind farm capacity to $\bm{u}^{\rm max} = \unit[(2.0, 3.0)]{p.u.}$ and truncate all generated samples of $\bm{\xi}\,|\,\bm{u}$ such that $\bm{\xi}\ge-\bm{u}$ and $\bm{\xi}\le\bm{u}^{\rm max} - \bm{u}$. 
The resulting collection of samples of $\bm{\zeta}$ with a single sample of $\bm{\xi}$ corresponds to what would be available in practice. We use 1500 of these generated samples for training and 500 for testing.

We define the following 5 cases that we use to demonstrate the proposed prescribed robust sets:
\begin{itemize}[left=0pt]
    \item \textit{Full}: Reference case for which the security interval of each wind farm covers the entire empirical forecast error support.
    \item \textit{90 Perc}: Reference case for which the security interval of each wind farm is fixed between the \unit[10]{\%} and \unit[90]{\%} percentile of the forecast error training data.
    \item \textit{Single}: Learning a single fixed uncertainty set without prescription, i.e., $\mathcal{M}^{\mu} = \bm{m}^{\mu}$ and $\mathcal{M}^{\sigma} = \bm{m}^{\sigma}$.
    \item \textit{P-All}: Learning $\bm{w}$ using all available forecast errors. As a result, the training is ignorant to the error distribution being conditional to $\bm{u}$.
    \item \textit{P-Cond}: Learning $\bm{w}$ assuming access to samples from the true conditional distribution. In each step $z$, given $\bm{u}^z$, we generated 200 new samples of $\bm{\zeta}$ using \cref{eq:forecast_error_sample_distribution}.
    \item \textit{P-Bins}: Learning $\bm{w}$ with each set of forecast error samples $\set{X}^z$ obtained by separating the available training samples of $\bm{u}$ into 10 bins of equal width, as in \cite{dvorkin2015uncertainty}, collecting the forecast error of each bin, and assigning each sample $\bm{u}^z$ the set of forecast error samples corresponding to its bin.
\end{itemize}

{Cases \textit{Full}, \textit{90 Perc.}, and \textit{Single} offer a comparison to methods with fixed robust sets.}
We set $v^{\rm max}=100$, although we typically observed satisfying convergence within 40 iterations. We used a mini-batch size of $z^{\rm max}=20$ and a learning rate of $\tau=10^{-6}$.
Finally, we set $\bm{w}^{\rm init}$ such that all entries of $\bm{M}^{\mu}$, $\bm{M}^{\sigma}$, and $\bm{m}^{\mu}$ are zero and the entries $\bm{m}^{\sigma}$ correspond to two times the empirical standard deviation of the training data.
All experiments are implemented in Python (see also Section~\ref{ssec:implementation}) and available online \cite{git_p-robust_dcopf}. We used a standard PC workstation with \unit[16]{GB} memory and an Intel i5 processor. 

\subsubsection{Cost-based box uncertainty set}
\label{ssec:case_study_coe_5bus}
We first train $\mathcal{M}_{\bm{w}}$ to optimize the expected cost of constraint exceedance (see Section~\ref{ssec:coe_solution_approach}).
We set $c_k^{\rm viol} = \unit[20k]{\nicefrac{\$}{MW}},\ \forall k$, which corresponds to the value of lost load estimated by New York ISO \cite{nyiso2019ancilary}. The average computation time for each epoch across all cases was around \unit[0.25]{s}, including sampling, solving the inner optimization problems, and computing the loss function and its gradients.

\begin{figure}
    \centering
    \includegraphics[width=0.95\linewidth]{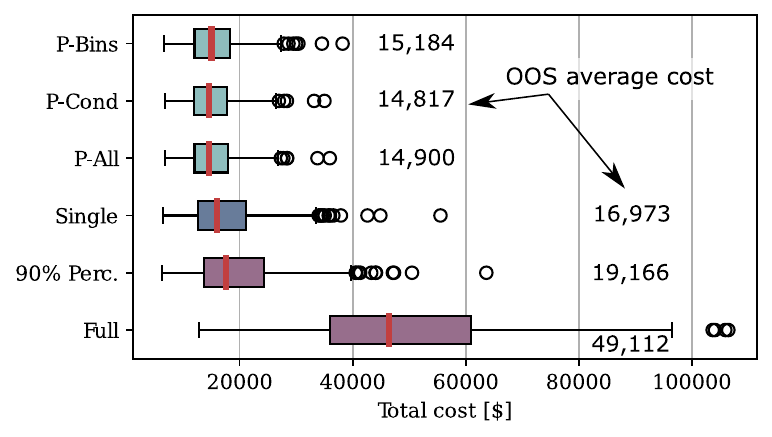}
    \caption{Out-of-sample (OOS) results for 5-bus system in configuration A. Optimized for expected cost of exceedance. The boxes extend from the first quartile to the third quartile. Red line shows the median. The whiskers extend from the box to the farthest data point within 1.5 times the inter-quartile range from the box.}
    \label{fig:5bus_cost_based_training}
\end{figure}

Fig.~\ref{fig:5bus_cost_based_training} shows the out-of-sample (OOS) average cost for the studied cases. 
The fully robust approach (\textit{Full}) is overly conservative, which we can infer from the fact that the lowest cost on this case are around the 25-percentile of all other cases. Also, it is often infeasible, leading to higher cost from infeasibility penalties (see Section~\ref{ssec:implementation} IR2). The fixed \textit{90 Perc.} case is less conservative and mostly avoids infeasibility, but is outperformed by the other cases. 
The \textit{Single} uncertainty set improves average cost relative to \textit{90 Perc}, but the set remains too small, as it tries to avoid infeasibility penalties.
Introducing the prescription step overcomes this problem. The prescriptive cases \textit{P-All}, \textit{P-Cond}, and \textit{P-Bins} further improve upon \textit{Single} by \unit[12.2]{\%}, \unit[12.7]{\%}, and \unit[10.5]{\%}, respectively. 
Case \textit{P-Cond} with access to the true conditional distribution slightly outperforms \textit{Single}. 
Case \textit{P-Bins} performs worse than \textit{P-Single}, which we attribute to the loss of correlation information in the binning of the forecast errors.  
Re-running the experiment without correlation ($\phi=0$) confirms this. Now, \textit{P-All}, \textit{P-Cond}, and \textit{P-Bins} result in average OOS cost of \unit[14,899]{\$}, \unit[14,851]{\$}, and \unit[14,981]{\$}, with \textit{P-Bins} now outperforming \textit{P-All} and only being slightly worse than the ideal \textit{P-Cond}.

\begin{figure}
    \centering
    \includegraphics[width=0.95\linewidth]{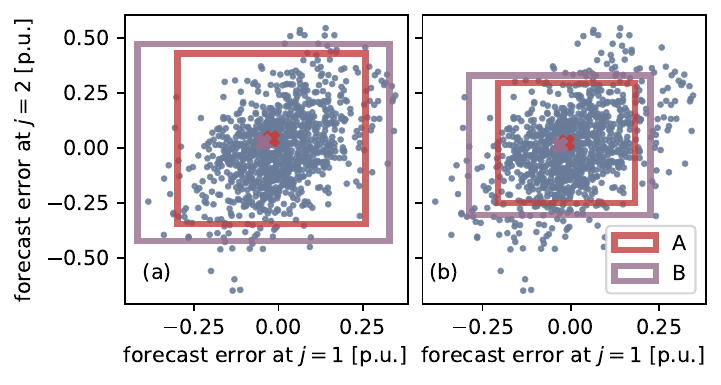}
    \caption{Prescribed sets from case \textit{P-All} (optimized for cost of exceedance) for (a) $\bm{\zeta}=(\bm{d}_0, \bm{u}_0)$ and (b) $\bm{\zeta}=0.7(\bm{d}_0, \bm{u}_0)$, each for the two network configurations shown in Fig.~\ref{fig:5bus_system}.}
    \label{fig:two_configs_with_prescription}
\end{figure}

Fig.~\ref{fig:two_configs_with_prescription} shows how a change of the system topology impacts the prescribed uncertainty set, highlighting the relevance of making the training problem-aware. By moving wind farm $j=2$ from bus 5 (configuration A) to bus 2 (configuration B), it can no longer rely on the direct balancing from the generator at bus 5. 
As a result, transmission lines are now more likely to exceed their remaining available margins in real time which leads to larger required safety intervals for both wind farms.

\subsubsection{Constraint-based box uncertainty set}
\label{ssec:case_study_poe_5bus}

\begin{figure}
    \centering
    \includegraphics[width=0.95\linewidth]{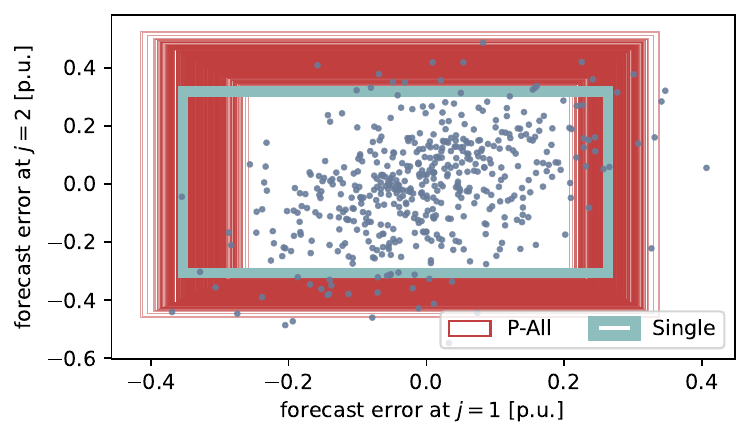}
    \caption{Uncertainty sets for cases \textit{P-All} and \textit{P-Single} optimized for probability of exceedance. For \textit{P-All} all prescriptions resulting from the test samples of $\bm{\zeta}$ are shown. (System configuration A.)}
    \label{fig:5bus_cvar_boxes}
\end{figure}

We now train $\mathcal{M}_{\bm{w}}$ such that the probability of constraint exceedance remains below \unit[99]{\%}, i.e., $\gamma=0.01$, as described in Section~\ref{ssec:poe_solution_approach}.
For this experiment we increase the learning rate to $\phi=10^{-5}$ and set $\lambda^{\rm init} = 100$ and $\kappa = 0.1$.
The time per iteration increased slightly to around \unit[0.67]{s} per epoch, which we can mainly attribute to the more complex loss function. 
For \textit{Single} the CVaR at convergence is 0.96 with an in-training probability of exceedance of 3.1\% and a testing probability of exceedance of 13\%. As discussed in Section~\ref{ssec:case_study_coe_5bus} above, this unsatisfying result can be explained with the algorithm avoiding infeasible uncertainty sets resulting in a set that is too small for many parametrizations (see Fig.~\ref{fig:5bus_cvar_boxes}).
For \textit{P-All}, on the other hand, the CVaR at convergence is 0.06 and we achieve an in-training probability of exceedance of 1\% and a testing probability of exceedance of 1.8\%. 
We note that the exact match of the in-training probability of exceedance with the target is a coincidence. The actual target of a CVaR equal to zero would lead to a smaller in-training probability of exceedance. However, the algorithm converges above this target as it again tries to avoid to create infeasible sets.
Fig.~\ref{fig:5bus_cvar_boxes} shows the resulting uncertainty sets for \textit{Single} and \textit{P-All} and highlights the improved ability of the set in the latter case to adapt to the current system state.
Comparing Figs.~\ref{fig:two_configs_with_prescription} and \ref{fig:5bus_cvar_boxes} we see that the uncertainty sets resulting from the two different performance targets are qualitatively similar.

\subsection{RTS 96-bus case}

We now test the suggested approach on more realistic data using the system provided by the Reliability Test System Grid Modernization Lab Consortium (RTS-GLMC \cite{rts_glmc_git}). 
This update of the RTS-96 test system has 73 buses, 120 transmission lines, 73 conventional generators, 4 wind farms, and 76 other resources (hydro, PV).
In our experiment we focus on the 4 wind farms as uncertain resources and treat hydro and PV injections as fixed negative demand, i.e. as part of $\bm{d}$.
The RTS-GLMC data set includes data for one year. To have a richer data set for the wind farms, we use the coordinates provided in the RTS-GLMC data set to map the 4 wind farms to the closest data points available in the extensive NREL WIND Toolkit \cite{draxl2015overview}.
We scale this data to fit the wind farms in the RTS-GLMC data set and obtain 7 years of wind power injections and realistic forecast errors. 
We select the data from 2012 to replace the wind data from the RTS-GLMC data set, as the yearly wind structure matches original data most closely (measured in terms of average deviation of hourly total wind injections). 
From the resulting 8760 available samples net-demand and wind-injection samples of the respective day-ahead data sets, we select 1500 for training and 500 for testing. 
We use all forecast errors for training as in \textit{P-All} above and focus on the analysis of the cost of exceedance-based training.

We select the same meta-parameters as for the 5-bus case, but reduce the mini-batch size to $z^{\rm max}=10$.
Training requires an average of \unit[28]{s} per epoch and we observe convergence after around 30 iterations. We note that around half of the time per iteration is spend on computing the gradient. This is expected because both larger parameter matrices $\bm{M}^{\sigma}$, $\bm{M}^{\mu}$ and inner optimization overproportionally increase the size of the computational graph from which the gradient is computed. 
However, because $\mathcal{M}_{\bm{w}}$ has to be trained only once offline, training time and resources are not a critical limiting factor.
{In addition, any relevant parameters that change more frequently than the time needed for re-training $\mathcal{M}_{\bm{w}}$, e.g., potential grid topology changes, can be included into context vector $\bm{\zeta}$.}

\begin{figure}
    \centering
    \includegraphics[width=0.95\linewidth]{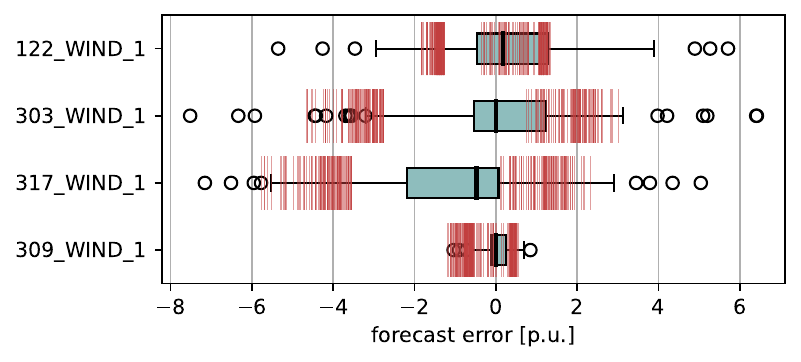}
    \caption{OOS results for the RTS 96 bus case using the \textit{P-All} training. Box plots show the distribution of the OOS forecast errors and red lines show the various security intervals obtained from the trained mapping. See Fig.~\ref{fig:5bus_cost_based_training} for box plot explanation.}
    \label{fig:rts96_coe_pall}
\end{figure}

In this case, the reference \textit{90 Perc.} leads to a expected out of sample cost of \unit[3,671,312]{\$} while the \textit{P-All} trained prescribed sets achieve a significant improvement of \unit[1,249,725]{\$}. (We performed an additional grid search to find a better percentile-based uncertainty set. The best result was attained for the 78-percentile with \unit[3,088,187]{\$}).
Fig.~\ref{fig:rts96_coe_pall} shows the OOS forecast errors alongside the collection of prescribed security intervals for the 4 wind farms in the system. 
We observe that the uncertainty sets are biased towards negative forecast errors. This can largely explained by the fact that the model chooses to curtail wind if the forecast is very high, which (i) leads to an insensitivity to upwards forecast errors and (ii) amplifies the fact that larger negative forecast errors are more likely at high wind power forecasts. 
This result highlights the advantage of internalizing the problem cost structure into the training.

\section{Conclusion}

This paper demonstrated an approach to compute uncertainty sets for robust optimal power flow that (i) are prescriptive, i.e., minimize the expected cost system cost, and (ii) are adaptive, i.e., are prescribed individually for each expected system state given by a vector of context parameters. 
Our approach to obtain these sets was inspired by \cite{wang2023learning} but additionally achieves property (ii).

The problem in \cref{eq:general_training_problem} opens a wide range of future research \cite{sadana2023survey}. For the approach studied in this paper we are pursuing the following avenues for further research. (i) Including richer context vectors. For example, in our data for the RTS96 case study, we observed a clear dependency of the forecast error distribution on the wind direction. Making the problem more context aware and better estimate conditional error distributions from real data should reveal impactful relations between system security and observable parameters, but is a non-trivial task \cite{sadana2023survey}. (ii) AC power flow. {Higher fidelity operational models, e.g., along the lines of \cite{lee2021robust,louca2016stochastic}, promise improved applicability in practice and improved system security.} (iii) Non-fixed recourse. Relaxing the second stage with a decision-making problem (or a tractable proxy, e.g., a trained neural network), should allow for a broader set of applications of the method. {(iv) Investigating scalability, e.g., through parallelization or differentiable optimization advances\mbox{, e.g., \cite{kotary2023folded}.}}

\section*{Acknowledgment}

The authors would like to thank Irina Wang and Bartolomeo Stellato of Pinceton University for their helpful discussions.

\bibliographystyle{IEEEtran}
\bibliography{literature}

\begin{thebibliography}{10}
\providecommand{\url}[1]{#1}
\csname url@samestyle\endcsname
\providecommand{\newblock}{\relax}
\providecommand{\bibinfo}[2]{#2}
\providecommand{\BIBentrySTDinterwordspacing}{\spaceskip=0pt\relax}
\providecommand{\BIBentryALTinterwordstretchfactor}{4}
\providecommand{\BIBentryALTinterwordspacing}{\spaceskip=\fontdimen2\font plus
\BIBentryALTinterwordstretchfactor\fontdimen3\font minus
  \fontdimen4\font\relax}
\providecommand{\BIBforeignlanguage}[2]{{%
\expandafter\ifx\csname l@#1\endcsname\relax
\typeout{** WARNING: IEEEtran.bst: No hyphenation pattern has been}%
\typeout{** loaded for the language `#1'. Using the pattern for}%
\typeout{** the default language instead.}%
\else
\language=\csname l@#1\endcsname
\fi
#2}}
\providecommand{\BIBdecl}{\relax}
\BIBdecl

\bibitem{roald2017chance}
L.~Roald and G.~Andersson, ``Chance-constrained ac optimal power flow,''
  \emph{IEEE Trans. Power Sys.}, vol.~33, no.~3, pp. 2906--2918, 2017.

\bibitem{lee2021robust}
D.~Lee \emph{et~al.}, ``Robust ac optimal power flow with robust convex
  restriction,'' \emph{IEEE Trans. Power Sys.}, vol.~36, no.~6, 2021.

\bibitem{bienstock2014chance}
D.~Bienstock \emph{et~al.}, ``Chance-constrained optimal power flow,''
  \emph{SIAM Review}, vol.~56, no.~3, 2014.

\bibitem{kazempour2018stochastic}
J.~Kazempour \emph{et~al.}, ``A stochastic market design with revenue adequacy
  and cost recovery by scenario: Benefits and costs,'' \emph{IEEE Trans. Power
  Sys.}, vol.~33, no.~4, pp. 3531--3545, 2018.

\bibitem{bertsimas2012adaptive}
D.~Bertsimas \emph{et~al.}, ``Adaptive robust optimization for the security
  constrained unit commitment problem,'' \emph{IEEE Trans. Power Sys.},
  vol.~28, no.~1, pp. 52--63, 2012.

\bibitem{liang2022inertia}
Z.~Liang \emph{et~al.}, ``Inertia pricing in stochastic electricity markets,''
  \emph{IEEE Trans. Power Sys.}, vol.~38, no.~3, 2022.

\bibitem{wang2023learning}
I.~Wang \emph{et~al.}, ``Learning for robust optimization,''
  \emph{arXiv:2305.19225}, 2023.

\bibitem{bertsimas2020predictive}
D.~Bertsimas and N.~Kallus, ``From predictive to prescriptive analytics,''
  \emph{Manag. Sci.}, vol.~66, no.~3, pp. 1025--1044, 2020.

\bibitem{lubin2015robust}
M.~Lubin \emph{et~al.}, ``A robust approach to chance constrained optimal power
  flow with renewable generation,'' \emph{IEEE Transactions on Power Systems},
  vol.~31, no.~5, pp. 3840--3849, 2015.

\bibitem{weinhold2023uncertainty}
R.~Weinhold and R.~Mieth, ``Uncertainty-aware capacity allocation in flow-based
  market coupling,'' \emph{IEEE Trans. Power Syst.}, 2023.

\bibitem{morales2023prescribing}
J.~M. Morales \emph{et~al.}, ``Prescribing net demand for two-stage electricity
  generation scheduling,'' \emph{Oper. Res. Perspect.}, vol.~10, p. 100268,
  2023.

\bibitem{bertsimas2018data}
D.~Bertsimas \emph{et~al.}, ``Data-driven robust optimization,'' \emph{Math.
  Program.}, vol. 167, pp. 235--292, 2018.

\bibitem{louca2018robust}
R.~Louca and E.~Bitar, ``Robust ac optimal power flow,'' \emph{IEEE Trans.
  Power Sys.}, vol.~34, no.~3, pp. 1669--1681, 2018.

\bibitem{louca2016stochastic}
------, ``Stochastic ac optimal power flow with affine recourse,'' in
  \emph{IEEE 55th Conference on Decision and Control (CDC)}.\hskip 1em plus
  0.5em minus 0.4em\relax IEEE, 2016, pp. 2431--2436.

\bibitem{golestaneh2018polyhedral}
F.~Golestaneh \emph{et~al.}, ``Polyhedral predictive regions for power system
  applications,'' \emph{IEEE Trans. Power Sys.}, vol.~34, no.~1, 2018.

\bibitem{lorca2016multistage}
A.~Lorca \emph{et~al.}, ``Multistage adaptive robust optimization for the unit
  commitment problem,'' \emph{Oper. Res.}, vol.~64, no.~1, pp. 32--51, 2016.

\bibitem{isuru2020piecewise}
M.~Isuru \emph{et~al.}, ``A piecewise-affine decision rule based stochastic ac
  optimal power flow approach,'' in \emph{Proc. of the IEEE PES Gen.
  Meeting}.\hskip 1em plus 0.5em minus 0.4em\relax IEEE, 2020.

\bibitem{morales2014electricity}
J.~M. Morales \emph{et~al.}, ``Electricity market clearing with improved
  scheduling of stochastic production,'' \emph{Eur. J. Oper. Res.}, vol. 235,
  no.~3, 2014.

\bibitem{vdvorkin2018setting}
V.~Dvorkin \emph{et~al.}, ``Setting reserve requirements to approximate the
  efficiency of the stochastic dispatch,'' \emph{IEEE Trans. Power Sys.},
  vol.~34, no.~2, pp. 1524--1536, 2018.

\bibitem{wang2013flexiramp}
B.~Wang and B.~F. Hobbs, ``Flexiramp market design for real-time operations,''
  in \emph{Proc. of the IEEE PES Gen. Meeting}.\hskip 1em plus 0.5em minus
  0.4em\relax IEEE, 2013.

\bibitem{garcia2021application}
J.~D. Garcia \emph{et~al.}, ``Application-driven learning: A closed-loop
  prediction and optimization approach applied to dynamic reserves and demand
  forecasting,'' \emph{arXiv:2102.13273}, 2021.

\bibitem{kotary2021end}
J.~Kotary \emph{et~al.}, ``End-to-end constrained optimization learning: A
  survey,'' \emph{arXiv:2103.16378}, 2021.

\bibitem{elmachtoub2022smart}
A.~N. Elmachtoub and P.~Grigas, ``Smart “predict, then optimize”,''
  \emph{Manag. Sci.}, vol.~68, no.~1, pp. 9--26, 2022.

\bibitem{agrawal2019differentiable}
A.~Agrawal \emph{et~al.}, ``Differentiable convex optimization layers,''
  \emph{Advances in neural information processing systems}, vol.~32, 2019.

\bibitem{donti2017task}
P.~Donti \emph{et~al.}, ``Task-based end-to-end model learning in stochastic
  optimization,'' \emph{Adv. Neural Inf. Process.}, vol.~30, 2017.

\bibitem{liang2022operation}
Z.~Liang \emph{et~al.}, ``Operation-adversarial scenario generation,''
  \emph{Electric Power Systems Research}, vol. 212, p. 108451, 2022.

\bibitem{wahdany2023more}
D.~Wahdany \emph{et~al.}, ``More than accuracy: end-to-end wind power
  forecasting that optimises the energy system,'' \emph{Electr. Power Syst.
  Res.}, vol. 221, p. 109384, 2023.

\bibitem{vdvorkin2023price}
V.~Dvorkin and F.~Fioretto, ``Price-aware deep learning for electricity
  markets,'' \emph{arXiv:2308.01436}, 2023.

\bibitem{bolognani2015fast}
S.~Bolognani and F.~D{\"o}rfler, ``Fast power system analysis via implicit
  linearization of the power flow manifold,'' in \emph{Proc. of the 53rd Annual
  Allerton Conference}.\hskip 1em plus 0.5em minus 0.4em\relax IEEE, 2015, pp.
  402--409.

\bibitem{bertsimas2011theory}
D.~Bertsimas \emph{et~al.}, ``Theory and applications of robust optimization,''
  \emph{SIAM review}, vol.~53, no.~3, pp. 464--501, 2011.

\bibitem{gorissen2013robust}
B.~L. Gorissen and D.~Den~Hertog, ``Robust counterparts of inequalities
  containing sums of maxima of linear functions,'' \emph{Eur. J. Oper. Res.},
  vol. 227, no.~1, pp. 30--43, 2013.

\bibitem{amari1993backpropagation}
S.-i. Amari, ``Backpropagation and stochastic gradient descent method,''
  \emph{Neurocomputing}, vol.~5, no. 4-5, pp. 185--196, 1993.

\bibitem{gower2019sgd}
R.~M. Gower \emph{et~al.}, ``Sgd: General analysis and improved rates,'' in
  \emph{Int. Conf. on Machine Learning}.\hskip 1em plus 0.5em minus 0.4em\relax
  PMLR, 2019, pp. 5200--5209.

\bibitem{dvorkin2015uncertainty}
Y.~Dvorkin \emph{et~al.}, ``Uncertainty sets for wind power generation,''
  \emph{IEEE Trans. Power Sys.}, vol.~31, no.~4, pp. 3326--3327, 2015.

\bibitem{mieth2023data}
R.~Mieth \emph{et~al.}, ``Data valuation from data-driven optimization,''
  \emph{arXiv:2305.01775}, 2023.

\bibitem{rockafellar2000optimization}
R.~T. Rockafellar \emph{et~al.}, ``Optimization of conditional value-at-risk,''
  \emph{Journal of Risk}, vol.~2, pp. 21--42, 2000.

\bibitem{zhang2022solving}
L.~Zhang \emph{et~al.}, ``Solving stochastic optimization with expectation
  constraints efficiently by a stochastic augmented lagrangian-type
  algorithm,'' \emph{INFORMS Journal on Computing}, vol.~34, no.~6, pp.
  2989--3006, 2022.

\bibitem{agrawal2019differentiating}
A.~Agrawal \emph{et~al.}, ``Differentiating through a cone program,''
  \emph{arXiv:1904.09043}, 2019.

\bibitem{domahidi2013ecos}
A.~Domahidi \emph{et~al.}, ``{ECOS}: {A}n {SOCP} solver for embedded systems,''
  in \emph{Proc. of the European Control Conference}, 2013, pp. 3071--3076.

\bibitem{wilder2019melding}
B.~Wilder \emph{et~al.}, ``Melding the data-decisions pipeline:
  Decision-focused learning for combinatorial optimization,'' in \emph{Proc. of
  the AAAI Conference on Artificial Intelligence}, vol.~33, no.~01, 2019.

\bibitem{matpowercase5}
\BIBentryALTinterwordspacing
{MATPOWER}. (2014) {CASE5 Power flow data}. [Online]. Available:
  \url{https://matpower.org/docs/ref/matpower5.0/case5.html}
\BIBentrySTDinterwordspacing

\bibitem{git_p-robust_dcopf}
\BIBentryALTinterwordspacing
R.~Mieth. (2023) {Prescribed Robust DC-OPF Code Supplement}. [Online].
  Available: \url{{https://github.com/mieth-robert/p-robust\_dcopf}}
\BIBentrySTDinterwordspacing

\bibitem{nyiso2019ancilary}
NYISO, ``Ancillary services shortage pricing,'' NYISO, Tech. Rep., 2019.

\bibitem{rts_glmc_git}
\BIBentryALTinterwordspacing
{Reliability Test System - Grid Modernization Lab Consortium}. [Online].
  Available: \url{github.com/GridMod/RTS-GMLC}
\BIBentrySTDinterwordspacing

\bibitem{draxl2015overview}
C.~Draxl \emph{et~al.}, ``Overview and meteorological validation of the wind
  integration national dataset toolkit,'' National Renewable Energy Lab.
  (NREL), Golden, CO (United States), Tech. Rep., 2015.

\bibitem{sadana2023survey}
U.~Sadana \emph{et~al.}, ``A survey of contextual optimization methods for
  decision making under uncertainty,'' \emph{arXiv:2306.10374}, 2023.

\bibitem{kotary2023folded}
J.~Kotary \emph{et~al.}, ``Folded optimization for end-to-end model-based
  learning,'' \emph{arXiv:2301.12047}, 2023.

\end{thebibliography}

\end{document}